\documentclass[a4paper,12pt]{amsart}
\usepackage{amssymb}
\author{Ivan V.~Arzhantsev}
\address{
Department of Higher Algebra\\
Faculty of Mechanics and Mathematics\\
Moscow State University\\
119992 Moscow, Russia}
\email{arjantse@mccme.ru}
\urladdr{http://mech.math.msu.su/department/algebra/staff/arzhan.htm}
\author{Natalia A. Tennova}
\address{
Department of Higher Algebra\\
Faculty of Mechanics and Mathematics\\
Moscow State University\\
119992 Moscow, Russia}
\email{tata@mccme.ru}
\thanks{Supported by RFBR grant MAC 03-01-06252, by CRDF
grant RM1-2543-MO-03 and by the RF President grant MK-1279.2004.1}

\title[On affinely closed homogeneous spaces]{On affinely closed
homogeneous spaces}
\date{June 24, 2004}
\dedicatory{To V.~N.~Latyshev on his 70-th anniversary}

\newcommand{\kk}{\Bbbk}

\newcommand{\Spec}{\mathop{\mathrm{Spec}}}

\newtheorem{theorem}{Theorem}
\newtheorem{corollary}{Corollary}

\newtheorem{lemma}{Lemma}

\theoremstyle{definition}
\newtheorem{definition}{Definition}

\theoremstyle{remark}
\newtheorem{remark}{Remark}

\begin{document}

\begin{abstract}
Affinely closed homogeneous spaces $G/H$, i.e.,
affine homogeneous spaces that admit only the trivial affine
embedding, are characterized for any affine algebraic
group~$G$. As a corollary, a description of affine $G$-algebras with finitely
generated invariant subalgebras is obtained.
\end{abstract}

\maketitle

\section{Introduction}
Let $G$ be an affine algebraic group over
an algebraically closed field $\kk$ of characteristic zero and $H$
an algebraic subgroup of $G$.  By Chevalley's Theorem, the
homogeneous space $G/H$ admits the canonical structure of a
quasiprojective variety. {\it An embedding} of the homogeneous
space $G/H$ is an algebraic $G$-variety $X$ with a base point
$x\in X$ such that the orbit $Gx$ is dense (and open) in $X$ and
the stabilizer $G_x$ equals $H$. We denote an embedding as
$G/H\hookrightarrow X$. We say that an embedding is {\it trivial}
if $Gx=X$. An embedding $G/H\hookrightarrow X$ is said to be {\it
affine} if the variety $X$ is affine. It is easy to show (for
example, see~\cite[Th.1.6]{pv}) that a space $G/H$ admits an
affine embedding if and only if $G/H$ is a quasi-affine variety,
or, equivalently, $H$ may be realized as the stabilizer of a
vector in a finite-dimensional $G$-module. In this case the
subgroup $H$ is said to be {\it observable} in $G$. An effective
description of observable subgroups in an affine algebraic group $G$ was
obtained by A.~Sukhanov~\cite{suh}.

\smallskip

The following definition was introduced in~\cite{at}.

\begin{definition}
A homogeneous space $G/H$ is called {\it affinely closed} if it admits
only the trivial affine embedding.
\end{definition}

An affinely closed homogeneous space is automatically affine. The
consideration of this class of homogeneous spaces is motivated by
the following question: when does "the stabilizer of a point $x$
on an affine $G$-variety $X$ equals $H$" imply "the orbit $Gx$ is
closed" ?

For a reductive group $G$ a homogeneous space $G/H$ is an affine variety
if and only if the subgroup $H$ is reductive (Matsushima's Criterion).
Note that for an arbitrary affine algebraic group
$G$ a description of affine homogeneous spaces $G/H$
in group-theoretic terms for the pair $(G,H)$ is an open problem,
for more information see~\cite[Ch.2]{gr}.

For reductive $G$ a description of affinely closed homogeneous spaces
follows directly from the result due to D.~Luna~\cite{lu}:

\begin{theorem}\label{tl}
Let $G$ be a reductive group. A homogeneous space
$G/H$ is affinely closed if and only if the subgroup $H$
is reductive and has a finite index in its normalizer $N_G(H)$.
Moreover, if $G$ acts on an affine variety  $X$ and the stabilizer of a point
$x\in X$ contains a reductive subgroup $H$ such that the group
$N_G(H)/H$ is finite, then the orbit $Gx$ is closed.
\end{theorem}

For example, for a maximal torus $T$ of a reductive group $G$ the
Weyl group $W=N_G(T)/T$ is finite, hence $G/T$ is affinely closed.
If $\rho:H\to SL(V)$ is an irreducible representation of a
semisimple group, then $SL(V)/\rho(H)$ is affinely closed (by the
Schur Lemma, the group $N_{SL(V)}(\rho(H))/\rho(H)$ is finite).
Thus the class of affinely closed homogeneous spaces is wide.

In~\cite{ar}, affinely closed homogeneous spaces of a reductive group $G$
play a key role in a classification of affine $G$-algebras such that
any invariant subalgebra is finitely generated.
Characterizations of complex affinely closed homogeneous spaces
of reductive groups in terms of compact transformation groups
and invariant algebras on compact homogeneous spaces are given in
\cite{lat} and~\cite{gl}.

The aim of this paper is to generalize the result of D.~Luna to the case
of an arbitrary affine algebraic group $G$ (Theorem~\ref{tm}) and to obtain
a classification of affine $G$-algebras with finitely generated invariant
subalgebras (Theorem~\ref{tmm}). Note that a characteristic-free variant
of these results for a solvable group $G$ is given in~\cite{ten}.


\section{A description of affinely closed spaces}
Let us fix a Levi decomposition $G=LG^u$ of the group $G$ in a
semidirect product of a reductive subgroup $L$ and the unipotent
radical $G^u$. By $\phi$ denote the homomorphism $\phi: G \to
G/G^u$. We shall identify the image of $\phi$ with $L$. Put
$K=\phi(H)$.

\begin{theorem} \label{tm}
The following conditions are equivalent:

\ {\rm (1)}\ $G/H$ is affinely closed;

\ {\rm (2)}\ $L/K$ is affinely closed.
\end{theorem}

\begin{proof}
The subgroup $H$ is observable in $G$
if and only if the subgroup $K$ is observable in $L$ \cite{suh},
\cite[Th.7.3]{gr}.

\smallskip

Suppose that $L/K$ admits a non-trivial affine embedding. Then
there are an $L$-module $V$ and a vector $v\in V$ such that the
stabilizer $L_v$ equals $K$ and the orbit boundary $Y=Z\setminus
Lv$, where $Z=\overline{Lv}$, is nonempty. Let $I(Y)$ be the ideal
in $\kk[Z]$ defining the subvariety $Y$. Recall that for an action
of an algebraic group $G$ on an affine variety $X$ any element
$f\in\kk[X]$ belongs to a finite-dimensional invariant subspace,
or, equivalently, $\kk[X]$ is a sum of its finite-dimensional
$G$-submodules. Thus there exists an $L$-submodule $V_1\subset
I(Y)$ that generates $I(Y)$ as an ideal. The inclusion
$V_1\subset\kk[Z]$ defines $L$-equivariant morphism $\psi:Z\to
V_1^*$ and $\psi^{-1}(0)=Y$. Then $L$-equivariant morphism $\xi:
Z\to V_2=V_1^*\oplus (V\otimes V_1^*)$, $z\to (\psi(z),
z\otimes\psi(z))$ maps $Y$ to the origin and is injective on the
open orbit in $Z$. Hence we obtain an embedding of $L/K$ in an
$L$-module such that the closure of the image of this embedding
contains the origin. Put $v_2=\xi(v)$. By the Hilbert-Mumford
Criterion, there is a one-parameter subgroup $\lambda: \kk^* \to
L$ such that $\lim_{t\to 0} \lambda(t)v_2=0$. Consider the weight
decomposition $v_2=v_2^{(i_1)}+\dots+v_2^{(i_s)}$ of the vector
$v_2$, where $\lambda(t)v_2^{(i_k)}=t^{i_k}v_2^{(i_k)}$. Here all
$i_k$ are positive.

By the identification $G/G^u=L$, one may consider $V_2$ as a
$G$-module. Let $W$ be a finite-dimensional $G$-module with a
vector $w$ whose stabilizer equals $H$. Replacing the pair $(W,
w)$ by the pair $(W\oplus (W\otimes W), w+w\otimes w)$, one may
suppose that the orbit $Gw$ intersects the line $\kk w$ only at
the point $w$. For a sufficiently large $N$ in the $G$-module
$W\otimes V_2^{\otimes N}$ one has $\lim_{t\to 0} \lambda(t)
(w\otimes v_2^{\otimes N})=0$ ($\lambda(\kk^*)$ may be considered
as a subgroup of $G$). On the other hand, the stabilizer of
$w\otimes v_2^{\otimes N}$ coincides with $H$. This implies that
the space $G/H$ is not affinely closed.

\smallskip

Conversely, suppose that $G/H$ admits a non-trivial affine
embedding. This embedding corresponds to a $G$-invariant
subalgebra $A\subset\kk[G/H]$ containing a non-trivial
$G$-invariant ideal $I$. Note that the algebra $\kk[L]$ may be
identified with the subalgebra in $\kk[G]$ of (left- or right-)
$G^u$-invariant functions, $\kk[G/H]$ is realized in $\kk[G]$ as
the subalgebra of right $H$-invariants, and $\kk[L/K]$ is the
subalgebra of left $G^u$-invariants in $\kk[G/H]$. Consider the
action of $G^u$ on the ideal $I$. By the Lie-Kolchin Theorem,
there is a non-zero $G^u$-invariant element in $I$. Thus the
subalgebra $A\cap\kk[L/K]$ contains the non-trivial $L$-invariant
ideal $I\cap\kk[L/K]$. If the space $L/K$ is affinely closed then
we get a contradiction with the following lemma.

\begin{lemma}\label{lv}
Let $L/K$ be an affinely closed space of a reductive group
$L$. Then any $L$-invariant subalgebra in $\kk[L/K]$ is finitely generated and
does not contain non-trivial $L$-invariant ideals.
\end{lemma}

\begin{proof}
Let $B\subset\kk[L/K]$ be a non-finitely generated
invariant subalgebra. For any chain $W_1\subset W_2\subset
W_3\subset\dots$ of finite-dimensional $L$-invariant submodules in
$\kk[L/K]$ with $\cup_{i=1}^{\infty} W_i=\kk[L/K]$, the chain of
subalgebras $B_1\subset B_2\subset B_3\subset\dots$ generated by
$W_i$ does not stabilize. Hence one may suppose that all
inclusions here are strict. Let $Z_i$ be the affine $L$-variety
corresponding the algebra $B_i$. The inclusion $B_i\subset
\kk[L/K]$ induces the dominant morphism $L/K\to Z_i$ and
Theorem~\ref{tl} implies that $Z_i=L/K_i$, $K\subset K_i$. But
$B_1\subset B_2\subset B_3\subset\dots$, and any $K_i$ is strictly
contained in $K_{i-1}$, a contradiction. This shows that $B$ is
finitely generated and, as proved above, $L$ acts transitively on
the affine variety $Z$ corresponding to $B$. But any non-trivial
$L$-invariant ideal in $B$ corresponds to a proper $L$-invariant
subvariety in $Z$.
\end{proof}

Theorem~\ref{tm} is proved.
\end{proof}

\begin{corollary}\label{c1}
Let $G/H$ be an affinely closed homogeneous space. Then for any affine
$G$-variety $X$ and a point $x\in X$ such that $Hx=x$, the orbit $Gx$
is closed.
\end{corollary}

\begin{proof}
The stabilizer $G_x$ is observable in $G$, hence $\phi(G_x)$ is
observable in $L$. The subgroup $\phi(G_x)$ contains $K=\phi(H)$,
and Theorems~\ref{tl} and~\ref{tm} imply that the space $L/\phi(G_x)$ is
affinely closed. By Theorem~\ref{tm}, the space $G/G_x$ is
affinely closed.
\end{proof}

 In particular, we get

\begin{corollary}\label{c2}
If $X$ is an affine $G$-variety and a point $x\in X$ is $T$-fixed,
where $T$ is a maximal torus of the group $G$, then the orbit $Gx$
is closed.
\end{corollary}


\section{$G$-algebras with finitely generated invariant
subalgebras}
Below {\it an affine algebra} over a field $\kk$ means
a finitely generated associative commutative $\kk$-algebra with
unit. Let $F$ be a subgroup of the automorphism group of an affine
algebra ${\mathcal A}$ and ${\rm rad}({\mathcal A})$ the set of nilpotent
elements of the algebra ${\mathcal A}$. Clearly, ${\rm rad}({\mathcal A})$
is an $F$-invariant ideal in ${\mathcal A}$.

\begin{lemma}\label{lr}
The following conditions are equivalent:

{\rm (1)}\ any $F$-invariant subalgebra in ${\mathcal A}$ is finitely generated;

{\rm (2)}\ any $F$-invariant subalgebra in ${\mathcal A}/{\rm rad}({\mathcal A})$
is finitely generated and $\dim {\rm rad}({\mathcal A})<\infty$.
\end{lemma}

\begin{proof}
Any finite-dimensional subspace in
${\rm rad}({\mathcal A})$ generates a finite-dimensional subalgebra in ${\mathcal A}$.
Hence if
$\dim {\rm rad}({\mathcal A})=\infty$, then the subalgebra generated by this
subspace is not finitely generated. On the other hand, the preimage in
${\mathcal A}$ of any non-finitely generated subalgebra in
${\mathcal A}/{\rm rad}({\mathcal A})$ is not finitely generated.

Conversely, suppose that (2) holds.
Then any subalgebra in ${\mathcal A}$ is generated by elements whose images
generate the image of this subalgebra in
${\mathcal A}/{\rm rad}({\mathcal A})$, and by a basis
of the radical of the subalgebra.
\end{proof}

By definition, a {\it $G$-algebra} is an affine algebra ${\mathcal A}$ with an
action (by automorphisms) of an algebraic group $G$ such that
any element $a\in {\mathcal A}$ is contained in a finite-dimensional
$G$-invariant subspace, where $G$ acts rationally.
Our aim is to describe all $G$-algebras such that any $G$-invariant
subalgebra is finitely generated.
By Lemma~\ref{lr}, we may assume that ${\rm rad}({\mathcal A})=0$.

Let $X=\Spec (\mathcal A)$ be the affine variety corresponding
to affine algebra ${\mathcal A}$ without nilpotents. To fix a
structure of $G$-algebra on ${\mathcal A}=\kk[X]$ is nothing else but to
fix an (algebraic) $G$-action on $X$. Define {\it the dimension}
$\dim{\mathcal A}$ of the algebra ${\mathcal A}$ as the dimension of the variety $X$.

\begin{lemma}\label{ld}
If $\dim{\mathcal A}\le 1$, then any subalgebra in ${\mathcal A}$ is finitely generated.
\end{lemma}

\begin{proof}
The case when $X$ is irreducible is considered
in~\cite[Prop.~2]{ar}. If $X=X_1\cup\dots\cup X_m$ is the decomposition on
irreducible components, then ${\mathcal A}$ is embedded into the direct sum of
the algebras $\kk[X_i]$, and any subalgebra in a summand is finitely generated.
Now it is easy to finish the proof by induction on $m$
considering the projection of ${\mathcal A}$ on $\kk[X_1]$.
\end{proof}

\begin{lemma}\label{lcom}
Suppose that $X=Z_1\cup Z_2$, where $Z_1$ and $Z_2$ are
closed invariant subvarieties. Then the following conditions are equivalent:

{\rm (1)}\ any invariant subalgebra in ${\mathcal A}$ is finitely generated;

{\rm (2)}\ any invariant subalgebra in $\kk[Z_1]$ and in $\kk[Z_2]$
is finitely generated.
\end{lemma}

\begin{proof}
If there is a non-finitely generated subalgebra in $\kk[Z_1]$,
then one may consider its preimage with respect to the restriction
homomorphism $\kk[X]\to\kk[Z_1]$. To prove the converse,
embed $\kk[X]$ in $\kk[Z_1]\oplus\kk[Z_2]$ and use the arguments
from the proof of the previous lemma.
\end{proof}

By Lemma~\ref{lcom}, one may assume that $G$ acts transitively
on the set of irreducible components of the variety $X$.

Below we generalize a construction from~\cite{ar} to the case
of non-connnected groups and reducible varieties.
Let $Y$ be a closed subvariety of $X$.
Consider a subalgebra
$$
 {\mathcal A}(X,Y)=\{ f\in\kk[X] \mid f(y_1)=f(y_2) \ \forall\, y_1,y_2\in Y \}.
$$

\begin{lemma}\label{lk}
If $Y$ contains an irreducible component of positive dimension that
does not coincide with any irreducible component of $X$, then
${\mathcal A}(X,Y)$ is not finitely generated.
\end{lemma}

\begin{proof}
Note that ${\mathcal A}(X,Y)=\kk\oplus I(Y)$. If ${\mathcal
A}(X,Y)$ is finitely generated, then one may assume that
generators $f_1,\dots,f_k$ are in $I(Y)$. Any monom in
$f_1,\dots,f_k$ of degree $s$ is in $I(Y)^s$. Hence it is
sufficient to prove that for some $l$ the space $I(Y)/I(Y)^l$ is
infinite-dimensional.

Let $Y=Y_1\cup\dots\cup Y_n$ and $X=X_1\cup\dots\cup X_m$ be the
decompositions on irreducible components, and $Y_1\subset X_1$,
$Y_1\ne X_1$, $\dim Y_1>0$. Suppose that $f\in I(Y)$ and $f$ is
not identically zero on $X_1$. Let ${\mathcal O}_{X_1,Y_1}$ be the
local ring of the subvariety $Y_1$ in $X_1$ and ${\mathcal I}$ its
maximal ideal. By the Nakayama Lemma, $\cap_{i=1}^{\infty} {\mathcal
I}^i=0$, hence after restriction to $X_1$ the element $f$ belongs
to ${\mathcal I}^{l-1}\setminus {\mathcal I}^l$ for some $l\ge 2$. Let $W$
be a subspace in $\kk[X]$ complementary to $I(Y_1)$. Note that
$\dim Y_1>0$ implies $\dim W=\infty$. The subspace $fW$ may be
considered as an infinite-dimensional subspace in ${\mathcal I}^{l-1}$
with zero intersection with ${\mathcal I}^l$. Hence $fW$ determines an
infinite-dimensional subspace in $I(Y)/I(Y)^l$.
\end{proof}

We conclude that any invariant subvariety $Y$ satisfying the conditions
of Lemma~\ref{lk}, determines the non-finitely generated invariant subalgebra
${\mathcal A}(X,Y)$ in $\kk[X]$.

\smallskip

By $G^0$ denote the connected component of unit of a group $G$.

\begin{theorem}\label{tmm}
Let ${\mathcal A}$ be a $G$-algebra without nilpotents with the non-trivial
induced action of the subgroup $G^u$. The following conditions are
equivalent:

\smallskip

{\rm (1)}\ any $G$-invariant subalgebra in ${\mathcal A}$ is finitely generated;

{\rm (2)}\ any $G$-invariant subalgebra in ${\mathcal A}$ does not contain
non-trivial $G$-invariant ideals;

{\rm (3)}\ any $L$-invariant subalgebra in ${\mathcal A}^{G^u}$ does not contain
non-trivial $L$-invariant ideals;

{\rm (4)}\ ${\mathcal A}=\kk[G/H]$, where $G/H$ is an affinely closed
homogeneous space;

{\rm (5)}\ ${\mathcal A}^{G^u}=\kk[L/K]$, where $L/K$ is an affinely closed
homogeneous space.
\end{theorem}

\begin{proof}
$(1)\Rightarrow (4)$ \ {\it Step 1.}\ Suppose that
the action $G:X$ is not transitive. The closure $Y$ of a $G$-orbit
on $X$ is an invariant subvariety and we may apply Lemma~\ref{lk}
with the only exceptions $Y=X$ and $\dim Y=0$. Hence (1) implies
that $G^0$ acts on any component $X_i$ either with an open orbit and the
boundary of this orbit is a finite set of points, or trivially. In these cases
the action $G^u:X$ is trivial~\cite[Th.3]{po}, see
also~\cite[Prop.4]{ar}.

{\it Step 2.}\ Suppose that the action $G:X$ is transitive.
Then $X=G/H$. If
$G/H$ admits a non-trivial affine embedding $G/H\hookrightarrow X'$,
then $\kk[X']$ is an invariant subalgebra in ${\mathcal A}$. By Step 1,
this subalgebra contains a non-finitely generated invariant subalgebra.

\smallskip

$(2)\Rightarrow (4)$ \ The absence of non-trivial invariant ideals in
${\mathcal A}$ implies that $X=G/H$, and the absence of non-trivial invariant ideals
in invariant subalgebras implies that $G/H$ does not admit non-trivial
embeddings.

\smallskip

Proofs $(4)\Rightarrow (1)$ and $(4)\Rightarrow (2)$ are analogous
to the proof of Lemma~\ref{lv} (one should use
Corollary~\ref{c1}). By the same arguments we get $(5)\Rightarrow
(3)$.

\smallskip

We know that $\kk[G/H]^{G^u}=\kk[L/K]$. Theorem~\ref{tm}
implies $(4)\Rightarrow (5)$.

\smallskip

$(3)\Rightarrow (2)$\ Let ${\mathcal B}$ be an invariant subalgebra in ${\mathcal A}$
and $I$ a non-trivial invariant ideal in ${\mathcal B}$. Then
$I\cap {\mathcal A}^{G^u}$ is a non-trivial invariant ideal in
${\mathcal B}\cap {\mathcal A}^{G^u}$.
This completes the proof of Theorem~\ref{tmm}.
\end{proof}

\begin{remark}
 1) \ The conditions of Theorem~\ref{tmm} are not
equivalent to the condition "any $L$-invariant subalgebra in ${\mathcal A}^{G^u}$
is finitely generated": one may consider $G=G^u=(\kk, +)$ acting
on $\kk[x,y]$ by the formula $(a, f(x,y)) \to f(x+ay,y)$.

\smallskip

\ 2)\ The implication $(1)\Rightarrow (5)$ is incorrect if
${\mathcal A}={\mathcal A}^{G^u}$, see Lemma~\ref{ld} and~\cite{ar}.

\smallskip

\ 3) \ The restriction ${\mathcal A}\ne{\mathcal A}^{G^u}$ is natural
because the case of reductive group actions was studied
in~\cite{ar} under the assumptions that $G$ is connected and $X$
is irreducible. But Lemma~\ref{lk} and above arguments show that
if a reductive $G$ acts transitively on the set of irreducible
components of $X$, then any invariant subalgebra in $\kk[X]$ is
finitely generated if and only if either the $G^0$-algebras
$\kk[X_i]$ for any irreducible component $X_i$
have (in the terminology of~\cite{ar}) type C or HV, or
$X=G/H$ and $G/H$ is affinely closed (in this case, the
$G^0$-algebras $\kk[X_i]$ may not have type N, see example 5)
below). Surprisingly, the restriction ${\mathcal A}\ne{\mathcal A}^{G^u}$
simplifies the main results.
\end{remark}

\begin{corollary}
Let $G/H$ be a quasi-affine homogeneous space. Suppose that $H$
does not contain $G^u$. Then either $G/H$ is affinely closed or there are
infinitely many pairwise nonisomorphic affine embeddings
$G/H\hookrightarrow X_i$ and a sequence of dominant equivariant morphisms:
$$
X_1\gets X_2\gets X_3\gets \dots \ \ .
$$
\end{corollary}

\begin{proof}
Let $G/H\hookrightarrow X$ be a non-trivial affine embedding
and $Y$ the complement to the open orbit in $X$.
The algerba ${\mathcal A}(X,Y)$ is not finitely generated.
Let $p_0\in I(Y)$ be a non-zero element and
$f_1,\dots,f_s$ a set of generators of $\kk[X]$. Put
$p_i=p_0f_i$, $i=1,\dots,s$ and extend the set $p_0,p_1,\dots,p_s$
to an (infinite) generating set
$p_0,p_1,\dots,p_s,h_1,h_2,\dots$ of the algebra ${\mathcal A}(X,Y)$.
The affine $G$-varieties $X_i$ corresponding to the algebras
$$
{\mathcal B}_i=\kk[<Gp_0,Gp_1,\dots,Gp_s,Gh_1,\dots,Gh_i>] \subset
{\mathcal A}(X,Y) \subset \kk[X] \subset \kk[G/H],
$$
\noindent define embeddings
$G/H\hookrightarrow X_i$ and the inclusions of algebras determine
the desired chain of dominant morphisms. In the sequence
$X_i$ there is a subsequence consisting of pairwise nonisomorphic
embeddings.
(By definition, an isomorphish of two embeddings sends the base point
to the base point and is the unique equivariant morphism
identical on the open orbit.)
\end{proof}


\section{Some remarks on affine embeddings of homogeneous spaces for
non-reductive groups}
In this section we consider elementary
examples that provide negative answers to some natural questions.

\smallskip

Let $V$ be a finite-dimensional $G$-module and $v\in V$.

\smallskip

1) {\it The orbit $Gv$ is closed, but the orbit $Lv$ is not closed.}\
Consider
$$
  V=\kk^2,\  v=(1,1), \
G=\left\{\left( \begin{array}{cc}
               t & a \\
               0 & 1
              \end{array} \right)\right\},\
L=\left\{\left( \begin{array}{cc}
               t & 0 \\
               0 & 1
              \end{array} \right)\right\}.
$$

\smallskip

2) {\it The orbit $Lv$ is closed, but the orbit $Gv$ is not
closed.} Consider $$
  V=\kk^2,\  v=(1,1), \
G=\left\{\left( \begin{array}{cc}
               t & a \\
               0 & t^{-1}
              \end{array} \right)\right\},\
L=\left\{\left( \begin{array}{cc}
               t & 0 \\
               0 & t^{-1}
              \end{array} \right)\right\}.
$$

\smallskip

3) {\it The space $G/H$ is affinely closed does not imply that
all $L$-orbits on $G/H$ are closed.}\ Consider
$$
G=\left\{\left( \begin{array}{cc}
               t & a \\
               0 & t^{-1}
              \end{array} \right)\right\},\,
H=L=\left\{\left( \begin{array}{cc}
               t & 0 \\
               0 & t^{-1}
              \end{array} \right)\right\},\,
x=\left( \begin{array}{cc}
          1 & 1 \\
          0 & 1
          \end{array} \right)
$$
$$
\Rightarrow \overline{LxH}\ne LxH.
$$

\smallskip

4) {\it If $G^0/H^0$ is affinely closed then
$G/H$ is affinely closed, but the converse is not true.}
\ One may take $G=SL(2)$ with any finite non-Abelian subgroup $H$.

\smallskip

5) {\it If $G^0/(H\cap G^0)$ is affinely closed then
$G/H$ is affinely closed, but the converse is not true.} \
Consider
$$
  G=N_{SL(2)}T, \ \ H=<
\left( \begin{array}{cc}
        0 & 1 \\
       -1 & 0
      \end{array} \right) >.
$$

\smallskip

6) {\it The condition $0\in\overline{Gv}$ does not imply that $0$
may be obtained as the limit of a one-parameter subgroup in $G$.}\
The corresponding example for a solvable
$G$ is given in~\cite[sec.11]{bir}.

\smallskip

An important open problem is to characterize affinely closed spaces
(for both reductive and non-reductive groups) over an algebraically
closed field of positive characteristic. In particular, it is not known
do we have here
Corollary~\ref{c1}. Some results in this direction
may be found in~\cite[sec.8]{ar}.


\end{document}